# ARITHMETICAL FUNCTIONS: INFINITE PRODUCTS

Garimella Rama Murthy,
Associate Professor,
International Institute of Information Technology,
Hyderabad---500032, AP, INDIA

**ABSTRACT**
In this technical report, certain interesting classification of arithmetical functions is proposed. The notion of "additively decomposable" and "multiplicatively decomposable" arithmetical functions is proposed. The concepts of arithmetical polynomials and arithmetical power series are introduced. Using these concepts, an interesting Theorem relating arithmetical power series and infinite products has been proved. Also arithmetical polynomials are related to probabilistic number theory. Furthermore some results related to the Waring problem are discussed.

1. **Introduction:**

    Homo-Sapien civilization exhibited considerable curiosity in working with natural numbers. The definition of a prime number naturally motivated mathematicians to investigate the nature of divisors of a composite number. They were able to identify "perfect numbers" as those composite numbers for which the sum of divisors is equal to the number itself. Thus, mathematicians introduced the first arithmetical function: sum of divisors of an integer [Dun]. With progress in mathematical research, after quite a few centuries, the following arithmetical functions were introduced and studied [Ram], [Kan] :

    (i)     d(n): Number of divisors of an integer 'n' ( including itself as a divisor )

    (ii)    $\sigma_s(n)$ : Sum of "s"th powers of divisors of an integer 'n'

    (iii)   $\pi(n)$ : Total number of prime numbers until 'n'

    (iv)    p(n) : Number of "partitions" of an integer 'n'

    (v)     $\varphi(n)$ : Total number of integers below 'n' that are relatively prime to 'n'.

    Using the fundamental theorem of arithmetic, Euler proposed the following identity: ( so called Euler product ) [Apo]

$$\prod_{p \epsilon P} \left[ 1 + \frac{1}{p^s} + \frac{1}{p^{2s}} + \cdots \right] = \sum_{n=1}^{\infty} \frac{1}{n^s} \ldots\ldots\ldots\ldots(1.1).$$

Also, in association with the partition function, Euler discovered the following identity for the generating function of the partition function [Apo] :

**Theorem:** For |x|< 1, we have

$$\prod_{m=1}^{\infty} \frac{1}{1-x^m} = \sum_{n=0}^{\infty} p(n)\, x^n \quad where\ p(0) = 1 \ldots\ldots\ldots\ldots(1.2)$$

Based on the efforts of earlier mathematicians, the author proposed a collection of novel arithmetical functions and identities in the spirit of equation (1.2). This technical report is a continuation of the efforts documented in [Rama1]. Also in this technical report, an effort is made to classify arithmetical functions and identify interesting families of arithmetical functions.

This technical report is organized as follows. In Section 2, we briefly provide an interesting classification of arithmetical functions. Then two novel number theoretic / mathematical concepts, namely
(i) Arithmetical Polynomials and (ii) Arithmetical Power Series are defined / introduced. Some interesting results related to these concepts are documented. Some interesting links to the ideas of probabilistic number theory are documented. In Section 3, based on the results in the journal paper [Rama1], some novel identities are proved. Unlike the arithmetical functions in [Rama1], the identities in this technical report deal with some of the well studied arithmetical functions such as $\sigma_s(n)$, $d(n), \omega(n)$. In Section 4, some future research directions are briefly outlined. In section 5, some results related to the Waring problem are documented.

**2. Classification of Arithmetical Functions: Innovative Number Theoretic Concepts:**

Number theorists are often concerned with sequences of real or complex numbers. In number theory such sequences are called *arithmetical functions*. Formally, we have the following definition [Apo].

**Definition:** A real or complex valued function defined on the positive integers is called an arithmetical function or a number-theoretic function.

An interesting observation related to the well known arithmetical functions leads to the following classification:

(I) **Memoryless Arithmetical Functions**: An arithmetical function, M(.) is "memoryless" if for each value of the independent variable, 'n', the function depends only on the unique factorization of 'n'. i.e. Let
$$n = p_1^{n_1} p_2^{n_2} \ldots \ldots p_l^{n_l}.$$
We have that
$$M(n) = f(p_1, p_2, \ldots, p_l\ ;\ n_1, n_2, \ldots, n_l\ )$$

i.e the function depends only on the primes occurring in the unique factorization of 'n' and their exponents. For instance d(n),$\omega(n)$ are memoryless arithmetical functions.

(II) **Arithmetical Functions with Memory**: An arithmetical function, L(.) has "memory" if for each value of the independent variable 'n', the function potentially depends on the past and/or future values of the independent variable 'n'. For example $\pi(n), \varphi(n)$ are arithmetical functions with memory.

In the following discussion, we consider "memoryless" arithmetical functions and provide an interesting classification of such functions.

(A) **Additively Decomposable Arithmetic Functions**: Given an integer 'n' with unique factorization
$$n = p_1^{n_1} p_2^{n_2} \ldots \ldots p_l^{n_l},$$
an arithmetical function, $\beta(.)$, is "Additively Decomposable" if and only if
$$\beta(n) = h(p_1, n_1) + h(p_2, n_2) + \cdots + h(p_l, n_l)$$
for some function "h(.)" for all 'n'.

(B) **Multiplicatively Decomposable Arithmetic Functions**: Given an integer 'n' with unique factorization
$$n = p_1^{n_1} p_2^{n_2} \ldots \ldots p_l^{n_l},$$
an arithmetical function, $\alpha(.)$, is "Multiplicatively Decomposable" if and only if
$$\alpha(n) = g(p_1, n_1) \, g(p_2, n_2) \ldots g(p_l, n_l)$$
for some function "g(.)" for all 'n'.

It is interesting to note that the following interrelationship exists between the above two classes of memoryless arithmetical functions. The following two claims follow readily.

**Claim 1:** Suppose $\alpha(.)$ is a multiplicatively decomposable arithmetical function. Then $log_e [ \alpha(.) ]$ ( composite function ) is an additively decomposable arithmetical function.

**Claim 2**: Suppose $\beta(.)$ is an additively decomposable arithmetical function. Then $\exp [\beta(.)]$ ( composite function ) is a multiplicatively decomposable arithmetical function.

Now, we want to understand the connection between concepts (A), (B) and the traditional concepts such as "multiplicative" arithmetical functions and "completely multiplicative" arithmetical functions. For the purposes of completeness, we provide the definitions below.

**Definition**: An arithmetical function f(.) is called "multiplicative" if f(.) is not identically zero and if
$$f( m \, n ) = f(m) \, f(n) \text{ whenever } (m,n) = 1$$
i.e. when integers 'm' and 'n' are relatively prime. Further f(.) is called "completely multiplicative " if we also have
$$f( m \, n ) = f(m) \, f(n) \text{ for all integers 'm' and 'n'}.$$

**Lemma 1:** A "memoryless" arithmetical function, f(.) is "multiplicatively decomposable" if and only if it is "multiplicative".

**Proof:** Let us consider two integers { S,T } which are relatively prime i.e. we have
$$S = p_1^{n_1} p_2^{n_2} \ldots \ldots p_l^{n_l}$$

$$T = q_1^{r_1} q_2^{r_2} \ldots \ldots q_j^{r_j}.$$

By definition of a memoryless multiplicatively decomposable arithmetical function i.e.
$$f(S) = g(p_1, n_1) \, g(p_2, n_2) \ldots g(p_l, n_l) \quad \text{and}$$

$$f(T) = g(q_1, r_1) \, g(q_2, r_2) \ldots g(q_j, r_j).$$

Since S and T are relatively prime, we have that
$$f(ST) = g(p_1, n_1) \, g(p_2, n_2) \ldots g(p_l, n_l) \, g(q_1, r_1) \, g(q_2, r_2) \ldots g(q_j, r_j)$$
$$= f(S) \, f(T)$$

Thus, f(.) is also multiplicative arithmetical function.

Now we would like to prove the other part i.e. Suppose f(.) is a memoryless multiplicative arithmetical function i.e. Given
$$N = p_1^{n_1} p_2^{n_2} \ldots \ldots p_l^{n_l} = M_1 M_2 \ldots M_l, \quad \text{it is clear that}$$
$\{M_i\}'s \text{ are relatively prime.}$ Thus by the definition, we have that
$$f(N) = f(M_1 M_2 \ldots M_l) = f(M_1) \, f(M_2) \ldots f(M_l)$$
$$= f(p_1, n_1) \, f(p_2, n_2) \ldots f(p_l, n_l) \quad .$$

Thus f(.) is multiplicatively divisible       Q.E.D.

**Remark 1:** It should be noted that the Euler totient function is an arithmetic function with memory that is also multiplicative. Thus, it is interesting to see if there is an arithmetical function with memory that is not multiplicatively divisible.

In view of claim 1, we propose the following classification of "additively decomposable arithemetical functions" based on those of multiplicative arithmetical functions.

**Claim 3:** Suppose f(.)_ is an additively decomposable arithmetical function. Then it is easy to see that if (M,N) =1 i.e. M, N are relatively prime, we have that
$$f(MN) = f(M) + f(N).$$

**Definition:** We say that f(.) is "**completely additively decomposable**" if we have
$$f(MN) = f(M) + f(N) \quad \text{for any two integers that are}$$
not necessarily relatively prime.

The following examples give such completely additively decomposable functions.

**Example 1**: Consider two integers M,N such that
$$M = p_1^{n_1} p_2^{n_2} \quad \text{and} \quad N = p_1^{s_1} q_2^{s_2}.$$
It is clear that M,N are not relatively prime.

Let f(n) be a memoryless arithmetical function which is defined as the sum of exponents of primes occurring in the unique factorization of "n". Thus
$$f(M) = n_1 + n_2 \quad \text{and} \quad f(N) = s_1 + s_2.$$
Also $\quad f(MN) = n_1 + n_2 + s_1 + s_2 = f(M) + f(N).$

Thus f(.) is completely additively decomposable.

**Example 2:** Let $M = p_1^{n_1} p_2^{n_2} \ldots \ldots p_l^{n_l}$ and the memoryless arithmetical function of interest be *"log (M )"*. It is easy to see that such an arithmetical function is completely additively decomposable.
- Claim 3 naturally motivates the following definition related to arbitrary arithmetical functions:

**Definition:** An arithmetical function f(.) is called "additive" if f(.) is not identically zero and if

$$f( m n ) = f(m) + f(n) \text{ whenever } (m,n) = 1$$

i.e. when integers 'm' and 'n' are relatively prime.
Further f(.) is called "completely additive" if we also have

$$f( m n ) = f(m) + f(n) \text{ for all integers 'm' and 'n'.}$$

Along the lines of Lemma 1, we have the following Lemma

**Lemma 2:** A "memoryless" arithmetical function, f(.) is "additively decomposable" if and only if it is "additive".

**Proof**: The proof follows similar reasoning as in Lemma 1 and is avoided for brevity.                                                                                                Q.E.D.
Along the lines of claim 1 and 2, we have the following claim relating arbitrary additive and multiplicative arithmetical functions:

**Claim 4:** Suppose $\alpha(.)$ is a multiplicative arithmetical function. Then $log_e [\, \alpha(.) \,]$ ( composite function ) is an additive arithmetical function. Also, Suppose $\beta(.)$ is an additive arithmetical function. Then $\exp[\, \beta(.) \,]$ ( composite function ) is a multiplicative arithmetical function.
In the following discussion, we introduce two novel number theoretic concepts based on arithmetical functions.

**Definition:** Given two arbitrary arithmetical functions, $\alpha(.) \; and \; \beta(.)$ which assume only integer values ( with $\beta(.)$ assuming positive integer values ) , the following polynomial
$$\sum_{n=0}^{M} \alpha(n) \, x^{\beta(n)} \quad \text{with } \alpha(0) = \beta(0) = 1 \text{ ( and x is a complex number )}$$
is defined as an **arithmetical polynomial**. Some special cases arise when $\alpha(n) = n$ and / or $\beta(n) = n$. Other special cases can also be studied.

**Definition :** Given two arbitrary arithmetical functions, $\alpha(.) \; and \; \beta(.) \; which \; assume$ only integer values ( with $\beta(.)$ assuming positive integer values ), the following power series
$$\sum_{n=0}^{\infty} \alpha(n) \, x^{\beta(n)} \quad \text{with } \alpha(0) = \beta(0) = 1 \text{ ( and x is a complex number )}$$
is defined as an **arithmetical power series.**
We expect that the study of this class of structured arithmetical functions will lead to new results. Specifically, when $\alpha(.) \; and \; \beta(.)$ correspond to some "interesting" arithmetical functions, we can study the properties of such polynomials. In this spirit let us consider the following example:

- **Arithmetical Polynomials/Power Series : Links to Probabilistic Number Theory:**

**Example 3:** Let $\omega(n)$ represent the number of distinct prime divisors of "n" ( occurring in the unique factorization of 'n' ). This arithmetical function was studied extensively in the field of probabilistic number theory. For instance Ramanujan studied this function in his famous paper with Hardy [Har] (the normal number of prime factors of a number n ).
Let us define the following arithmetical polynomial i.e.
$$J(x) = 1 + \sum_{n=1}^{M} x^{\omega(n)}.$$

- Let "L" be the maximum number of prime divisors of any integer less than or equal to "M". Thus by using the definition of J(x), we have the following equivalent representation:

$J(x) =$ 1+ ( # of integers less than or equal to M with one prime divisor ) $x$
  + ( # of integers less than or equal to M with two prime divisors) $x^2$
  +……………………………………..
  +(# of integers less than or equal to M with 'L' prime divisors) $x^L$.

From the above, it is clear that $J(x)|_{x=1} = (M + 1)$. Thus, J(x) − (M+1) has a zero at x=1. It is interesting to see if all the zeroes are real for such a sequence of polynomials

- Now, we normalize the above polynomial with (M+1) :
$$\frac{J(x)}{(M+1)} = q_0 + q_1 x + q_2 x^2 + \cdots + q_L x^L = H(x),$$
where $\{q_i's\}$ are probabilities with an interesting interpretation.

- Let us provide the probabilistic interpretation of the above polynomial. Let us define a discrete random variable which assumes "L" values and let the probability mass function assume values in the set $\{q_0, q_1, q_2, \ldots\ldots, q_L\}$. It is clear that H(x) is the moment function of such discrete random variable [Pap].
- Thus, we have a sequence of probability distribution functions and the associated moment functions. It is interesting to investigate the properties of such functions.
More generally, we now consider the case of an arithmetic polynomial motivated by the above example and provide the links to the research area of "probabilistic number theory".

**Example 4:** Consider the arithmetic polynomial $1 + \sum_{n=1}^{M} x^{\beta(n)}$, where $\beta(.)$ is an *arbitrary* arithmetic polynomial. Let $\beta(.)$ assume values in the set $\{s_1, s_2, \ldots s_L\}$ such that $s_1 < s_2 < \cdots < s_L$. Also, let the exponent of $'x'$ be such that the value $s_j$ is assumed $t_j$ times for $1 \leq j \leq L$. Thus, the arithmetic polynomial is of the following form:
$$1 + t_1 x^{s_1} + t_2 x^{s_2} + \cdots + t_L x^{s_L}$$
such that $\sum_{i=1}^{L} t_i = M$. Now, on normalization by (M+1), we have

$$G(x) = \frac{1}{(M+1)} + \frac{t_1}{(M+1)} x^{S_1} + \cdots\cdots + \frac{t_L}{(M+1)} x^{S_L} = \frac{1}{(M+1)} + \sum_{i=1}^{L} q_i\, x^{S_i}.$$

- Thus, as in the above example G(x) represents the moment function of an associated discrete random variable [Pap].

By abstracting the logic behind the above two examples, we provide a general example.

**Example 5:** Let us consider, the arithmetical power series $\sum_{n=0}^{\infty} \alpha(n)\, x^{\beta(n)}$,

with $\alpha(0) = 1$ and $\alpha(.)$ assuming positive values ( not necessarily integers ).
Furthermore let $\sum_{n=0}^{\infty} \alpha(n) < \infty$. Under these reasonable assumptions, the arithmetical power series reduces to a power series of the following form:
$\sum_{n=0}^{\infty} a_n\, x^n$ with the assumption $\sum_{n=0}^{\infty} a(n) = S < \infty$.
As in the above two examples, normalize the numerical power series with $S$. Thus, we have the following power series:
$$K(x) = \frac{a_0}{S} + \sum_{n=1}^{\infty} \frac{a_n}{S}\, x^n = \sum_{n=0}^{\infty} q_n\, x^n$$
where $q_n's\ are$ probabilities. Thus K(x) represents the moment function of the associated discrete random variable.

- We expect that these examples lead to detailed results in probabilistic number theory.

In the journal paper [Rama1], various novel arithmetical functions are defined with the hope that they will be studied in future. In that work arithmetical power series based on these novel arithmetical functions are related to some infinite products. In the next section we consider arithmetical power series based on well studied arithmetical functions and relate them to associated infinite products.

3. **Well Studied Arithmetical Functions : Infinite Products**:

In this section, utilizing the Fundamental Theorem of Arithmetic (i.e the unique factorization of any integer using the prime numbers), several identities are derived. It is clear that the identities hold true from the point of view of SYMBOLIC ALGEBRA. Using standard arguments from analysis, it can easily be shown that *most of the identities hold true from the point of view of convergence.* Detailed analytic arguments are avoided for the sake of brevity. Lemma A below is repeated from [Rama1].

**LEMMA A:** Consider the following infinite product in unknown variable 'x' ( on the left hand side of the following identity ). It is true that

$$\prod_{p \in P} \left(1 + \frac{x}{p^k} + \frac{x}{p^{2k}} + \ldots\right) = 1 + \sum_{n=1}^{\infty} \frac{x^{\varpi(n)}}{n^k}, \quad \ldots\ldots\ldots\ldots\ldots(3.1)$$

where P is the set of prime numbers and $\varpi(n)$ is the number of distinct prime divisors of "n".

**Proof:**

The Left Hand Side (LHS) of the above equation, (3.1) becomes

$$\text{L.H.S} = \prod_{p \in P} \left(1 + \frac{x/p^k}{1 - \frac{1}{p^k}}\right)$$

Using the unique factorization of any integer, it is easy to see that the above equality holds true from the point of view of symbolic algebra. Now we consider its truth from the convergence view point.

**Typical Convergence Argument:**

Consider the power series on the right hand of equation (1). Upper bound this power series by the following power series i.e.

$$\left(1 + \sum_{n=1}^{\infty} \frac{x^{\varpi(n)}}{n^k}\right) < \left(1 + \sum_{n=1}^{\infty} x^n\right) \ldots\ldots\ldots(1.5)$$

It is easy to see that the power series defined on RHS of (1.5) converges absolutely on the interval (-1,1). Thus, invoking standard Theorem from real analysis, it can be seen that the power series on the left hand side of (1) converges. Q.E.D.

- Now we attempt to derive similar identities, where the arithmetic power series on the right hand side of the identity involves some well studied arithmetical functions. For instance, suppose the unique factorization of an arbitrary integer N is given as follows:

$$N = p_1^{n_1} p_2^{n_2} \ldots \ldots p_l^{n_l}.$$

It is well known that the sum of $t^{th}$ powers of all divisors of N ( including N itself) is given by

$$\sigma_t(N) = \left(1 + p_1^t + p_1^{2t} + \cdots + p_1^{n_1 t}\right) \ldots \ldots \left(1 + p_l^t + p_l^{2t} + \cdots + p_l^{n_l t}\right).$$

Using the above equation, the following identity is derived.

**LEMMA B:** Consider the following infinite product in unknown variable 'x' ( on the left hand side of the following identity ). It is true that

$$\prod_{p \, \varepsilon \, P}\left[1 + \frac{(1+p^t)x}{p^s} + \frac{(1+p^t+p^{2t})x}{p^{2s}} + \cdots + \frac{(1+p^t+\cdots+p^{n_1 t})x}{p^{n_1 s}} + \cdots\right] = 1 + \sum_{n=2}^{\infty} \frac{\sigma_t(N) \, x^{\omega(N)}}{N^s}$$

where P is the set of prime numbers and $\varpi(n)$ is the number of distinct prime divisors of "n".

**Proof:** It is easy to see that the most general term on the right hand side of the identity is of the form:

$$\frac{\left(1 + p_1^t + \cdots + p_1^{n_1 t}\right)\left(1 + p_2^t + \cdots + p_2^{n_2 t}\right) \ldots \left(1 + p_l^t + \cdots + p_l^{n_l t}\right) \, x^{\omega(N)}}{N^s}$$

Thus, the identity readily follows using the definition of $\sigma_t(N)$ and $\omega(N)$.
Q.E.D.

- Using a similar reasoning as in the above Lemma, the following identity also follows.

**Lemma C:** Consider the following infinite product in unknown variable 'x' ( on the left hand side of the following identity ). It is true that

$$\prod_{p \, \varepsilon \, P}\left[1 + \frac{(1+1)x}{p^s} + \frac{(1+2)x}{p^{2s}} + \cdots + \frac{(1+n_1)x}{p^{n_1 s}} + \cdots \right] = 1 + \sum_{n=2}^{\infty} \frac{d(N) \, x^{\omega(N)}}{N^s}$$

where P is the set of prime numbers, $\varpi(N)$ is the number of distinct prime divisors of "N" and d(N) is the number of divisors of N.

**Proof:** Avoided for brevity.

- Once again, the following Lemma also follows based on the unique factorization of "N".

**Lemma D:** Consider the following infinite product in unknown variable 'x' ( on the left hand side of the following identity ). It is true that

$$\prod_{p \, \varepsilon \, P}\left[1 + \frac{x}{p^s} + \frac{x^{2^t}}{p^{2s}} + \cdots + \frac{x^{n_1^t}}{p^{n_1 s}} + \cdots \right] = 1 + \sum_{n=2}^{\infty} \frac{x^{L(N)}}{N^s}$$

where P is the set of prime numbers, $L(N)$ is the sum of "t"th powers of exponents appearing in the unique factorization of "N" i.e.
$$L(N) = n_1^t + n_2^t + \cdots + n_l^t.$$

For each prime "p", the infinite sum on the LHS of the above identity is related to a theta function. Particularly when "t" equals to two, the connection with theta functions is apparent.

- It is seen that all the above identities (leading to novel arithmetical functions) are effectively derived based on the unique factorization in the integer ring. By abstracting the essential idea (in all the above identities), the following most general identity is derived.

**THEOREM E:** Consider the following infinite product

$$\prod_{p \in P}(1 + \frac{\theta(p,1) \, x^{K(p,1)}}{p^k} + \frac{\theta(p,2) \, x^{K(p,2)}}{p^{2k}} + \frac{\theta(p,3) x^{K(p,3)}}{p^{3k}} \cdots)$$

$$= 1 + \sum_{n=2}^{\infty}(\theta(p_1,n_1)\ldots\theta(p_l,n_l)) \, x^{K(p_1,n_1)+\ldots+K(p_l,n_l)})/n^k,$$

$$= 1 + \sum_{n=2}^{\infty} \frac{\alpha(n) \, x^{\beta(n)}}{n^k},$$

where the integer "n" has the following unique factorization
$$n = p_1^{n_1} \, p_2^{n_2} \ldots p_l^{n_l} \quad \text{and}$$

$$\alpha(n) = \theta(p_1, n_1)\, \theta(p_2, n_2) \ldots \theta(p_l, n_l),$$

$$\beta(n) = K(p_1, n_1) + K(p_2, n_2) + \cdots + K(p_l, n_l)$$

**Proof:** As in the previous identities, by considering an arbitrary term on LHS, It can be seen that the identity holds true.

- From the above most general identity, it is seen that $\alpha(.)$ is a multiplicatively decomposable memoryless arithmetical function, whereas $\beta(.)$ is an additively decomposable memoryless arithmetical function. In view of claims 1 and 2, these types of memoryless arithmetical functions are related to one another.

**REMARK:** Utilizing the logical procedure adopted in this identity (and thus all previous identities), various new arithmetical functions can be defined and associated with infinite products (in variable x). Detailed enumeration of all possible identities is avoided for brevity.

In the following section, some future research directions are briefly outlined.

### 4. Future Research Directions:

- Riemann zeta function was proposed by allowing the independent variable, 'k' in Euler sums i.e. $\sum_{n=1}^{\infty} \frac{1}{n^k}$ assume "complex" ( complex numbers ) values. Riemann hypothesis was proposed as a conjecture dealing with the zeroes of the Zeta function i.e. $\zeta(s)$.

In all the identities proposed in [Rama], the independent variable, 'x' can be allowed to assume complex values. Infact, for instance, we have the following identity

$$\prod_{p \in P} \left(1 + \frac{x}{p^s} + \frac{x}{p^{2s}} + \ldots \right) = 1 + \sum_{n=1}^{\infty} \frac{x^{\varpi(n)}}{n^s}.$$

Let both the independent variables { x, s } be allowed to assume complex values. Thus, we have the following function of two complex variables:

i.e. $1 + \sum_{n=1}^{\infty} \frac{x^{\varpi(n)}}{n^s} = R(x, s)$. It is easily seen that

$$R(1, s) = \zeta(s).$$

Thus, the nature of zeroes of such a function ( of two complex valued variables ) could have interesting structure like those of the Riemann zeta function.

**Note:** All the other identities derived in [Rama1], and this technical report lead to interesting functions of two complex variables. It is interesting to see if all the transcendental functions satisfy interesting functional equations ( as that in the case of the Riemann zeta function ).

- In the context of digital signal processing, most of the arithmetical functions constitute "digital" signals. Thus the results of digital signal processing could be utilized to derive new results related to the arithemetical functions.

- Ramanujan derived congruence properties of the "partition" function [Kan]. We expect interesting congruence properties for the arithmetical functions proposed in [Rama1] and this technical report. Also, detailed arithmetical information about these arithemetical functions will be extracted hopefully.

5. **On the Waring Problem**: **Generating Functions**:

    Bachet conjectured that every positive integer can be written as the sum of squares of atmost 4 integers. Lagrange provided the first proof of such a Theorem. This theorem naturally leads to the following arithmetical function:

    J(n) : Number of ways in which an integer 'n' can be expressed as the sum of squares of four integers.

    Our goal is to provide an expression for the generating function associated with such an arithmetical function. In this connection, let us consider the following theta function:

$$K(x) = \sum_{n=-\infty}^{\infty} x^{n^2}$$

Theta functions are very well studied by mathematicians such as Jacobi, Ramanujan etc. In fact Jacobi proved the elegant four square theorem comprehensively summarizing the arithmetical information about J(n).

**Goal:** We are interested in arriving at algebraic and analytic information about J(n). We have the following interesting algebraic result.
.
$$[K(x)]^4 = \sum_{p=-\infty}^{\infty} \sum_{q=-\infty}^{\infty} \sum_{r=-\infty}^{\infty} \sum_{s=-\infty}^{\infty} x^{p^2+q^2+r^2+s^2} = 1 + \sum_{m=1}^{\infty} J(m) \, x^m = I(x)$$

Thus, the power series above provides the generating function, I(x) for the arithmetical function { J(.) }.

(A) Euler first provided a generating function for the partition function, p(n). Using the results of analysis, the asymptotic order of p(n) was derived by Rademacher improving results of Ramanujan and Hardy. Also, the so called Hardy- Littlewood circle method is related to the Cauchy residue theorem.

- Using the well known approaches of analysis associated with the generating function of arithmetical functions, we hope to derive new results for the arithmetical function { J (n) } [Rama2].

(B) It is well known from the theory of infinite series that the following holds true

$$\left[ \sum_{n=0}^{\infty} a_n\, x^n \right]^2 = \sum_{n=0}^{\infty} ( a_0\, a_n + a_1\, a_{n-1} + \cdots + a_n\, a_0 )\, x^n$$

In the language of digital signal processing, the function { r(n) } constitutes the correlation function of { a(n) }. Ramanujan, in his classic paper [Ram], studied an arithmetical function related to the correlation function of $\sigma_s(n)$. Thus, we are naturally motivated to study the correlation function of { J(n) } denoted by { R(n) }. It is thus clear that the generating function for the correlation function of { J(n) } is given by considering the eighth power of K(x) ( or equivalently the square of I(x) ).. Hence we have that

$$[\, I(x)\, ]^2 = \sum_{n=0}^{\infty} R_n\, x^n$$

In this connection, we realize that Fermat first considered the representation of primes as sum of squares. He showed that every prime of the form 4 n + 1 can be represented as a sum of two integral squares in a single manner. He also reasoned that a prime of the form 4n+3 can never be written as sum of squares of two integers. Thus, we are naturally led to the following arithmetical function: { S(n) }

S(n) : Number of ways in which 'n' can be represented as the sum of squares of two integers.

In view of Fermat's result, S(n) is equal to 'zero' for primes of the form 4n+3 and is equal to 'one' for primes of the form 4n+1. The generating function for { S(n) } is given by the square of K(x). We are thus naturally led to the conclusion that "Correlation function of S(n) equals the arithmetical function J(n) whose arithmetical properties are thoroughly provided by the four square theorem of Jacobi. It may be possible to derive arithmetical properties of S(n) from those of J(n).

Based on the arithmetical function { S(n) }, we define the following novel arithmetical function

t(n) : Number of primes of the form 4n+1 below "n".

As in the case of Prime Number Theorem, we would like to determine the asymptotic order of t(n). Considering a special case of Lambert series, Wiener proved the prime number theorem. The result uses the ideas related to proof of Tauberian theorems. We expect that similar approach enables determining the asymptotic order of t(n) [Kno].

- Waring considered a generalized version of the above problem. He conjectured that every number can be represented as the sum of a limited number of cubes, fourth or higher powers of integers. That such representations exists was proved by Hilbert. In connection with the study of Waring's problem mathematicians introduced the following generalization of the theta function.

$$C^{(s)}(x) = \sum_{n=-\infty}^{\infty} x^{n^s}$$

In the above, we first consider the case where "s" is an even integer. The case where "s" is an odd integer is considered in [Rama2].

Generalizing the above discussion related to the arithmetical function { J(n) }, we consider higher powers of $C^{(s)}(x)$ i.e.

$$[\, C^{(s)}(x) \,]^t$$

The higher powers of $C^{(s)}(x)$ naturally lead to the introduction of the following very general arithmetical function associated with the Waring's problem.

$W_{(t)}^{(s)}(m)$ : Number of ways in which an integer "m" can be written

as the "s"th power of "t" integers i.e.

$$m = n_1^s + n_2^s + \cdots + n_t^s.$$

Let "B" be the set of integers expressible as the s"th power ( where "s' is an even integer ) of "t" integers.

- The generating function of such a general arithmetical function is given by
$$U_{(t)}^{(s)}(x) = U(x;t;s) = [C^{(s)}(x)]^t = \sum_{m \in B} W_{(t)}^{(s)}(m) x^m,$$

**Goal:** To thoroughly investigate the arithmetic, algebraic and analytic properties of the above general arithmetical function associated with the Waring problem.. In this connection, the following functional equation naturally follows:

**Lemma G:**   U ( x ; t+r ; s ) = U ( x; t; s ) U ( x; r; s )

**Proof:** It follows based on the fact that the generating function of $W_{(t)}^{(s)}(m)$ is the 't'th power of generalized theta function. Details are avoided for brevity   Q. E. D.

**Corollary:**
$$W_{(t+r)}^{(s)}(m) = W_{(t)}^{(s)}(m) * W_{(r)}^{(s)}(m),$$

where "*" denotes the convolution of infinite sequences associated with the general arithmetical function ( naturally connected with the Waring problem ).

- **Relevance of Infinite Series Related Results** :

  The essential idea in the above discussion is to consider the higher powers of theta function and the generalized theta function. It is thus natural to invoke results available from the theory of infinite series ( dealing with higher powers of a given infinite series ). From [Kno], it is well known that

$$\left[\sum_{n=0}^{\infty} a_n x^n\right]^k = \sum_{n=0}^{\infty} a_n^{(k)} x^n,$$

where the coefficients $a_n^{(k)}$ are constructed from the coefficients $a_n$ in a perfectly determinate manner. Recurrence formula for the evolution of $a_n^{(k)}$ are already derived. Thus the hope is that the arithmetical functions $W_{(t)}^{(s)}(m)$ ( for fixed 's' and increasing values of 't' ) satisfy a linear / non-linear recurrence equation.

**Goal:** Using the results well known for difference equations and an initial condition, it could be possible to extract arithemetical information connected with the arithmetical function $W_{(t)}^{(s)}(m)$. The author is actively pursuing these ideas [Rama2].

6. **Conclusions:**

In this technical report, classification of arithmetical functions is attempted. Additively and multiplicatively decomposable arithmetical functions are defined. The novel concepts of arithmetical polynomials and arithmetical power series are introduced. An interesting Theorem associated with arithmetical power series and infinite products is proved. Arithmetical power series are linked to the research area of probabilistic number theory. It is expected that this technical report leads to detailed investigation of "additively decomposable" memoryless arithmetical functions. Also some results related to the Waring problem are provided.